\documentclass[11pt,english]{article}
\usepackage[T1]{fontenc}
\usepackage[latin9]{inputenc}
\usepackage[a4paper]{geometry}
\usepackage{babel}
\usepackage{amsmath}
\usepackage{amssymb}
\usepackage[unicode=true,pdfusetitle,
 bookmarks=true,bookmarksnumbered=false,bookmarksopen=false,
 breaklinks=false,pdfborder={0 0 1},backref=section,colorlinks=false]
 {hyperref}
\usepackage{breakurl}

\makeatletter

\newcommand{\noun}[1]{\textsc{#1}}

\usepackage{babel}

\makeatother

\begin{document}

\title{Literature survey on low rank approximation of matrices\thanks{Supported by National Board of Higher Mathematics, India.}}

\author{N. Kishore Kumar\thanks{BITS-Pilani Hyderabad Campus; India, Email: thimmaki@gmail.com,naraparaju@hyderabad.bits-pilani.ac.in}$^{\star}$
and J. Schneider\thanks{Friedrich Schiller Universität Jena, Germany; Email: jan.schneider@uni-jena.de }$^{\dagger}$
\\
 }

\date{~}

\maketitle
 
\begin{abstract}
Low rank approximation of matrices has been well studied in literature.
Singular value decomposition, QR decomposition with column pivoting,
rank revealing QR factorization (RRQR), Interpolative decomposition
etc are classical deterministic algorithms for low rank approximation.
But these techniques are very expensive ($O(n^{3})$ operations are
required for $n\times n$ matrices). There are several randomized
algorithms available in the literature which are not so expensive
as the classical techniques (but the complexity is not linear in $n$).
So, it is very expensive to construct the low rank approximation of
a matrix if the dimension of the matrix is very large. There are alternative
techniques like Cross/Skeleton approximation which gives the low-rank
approximation with linear complexity in $n$. In this article we review
low rank approximation techniques briefly and give extensive references
of many techniques. \medskip{}
\medskip{}
\end{abstract}

\section{Introduction}

The low rank matrix approximation is approximating a matrix by one
whose rank is less than that of the original matrix. The goal of this
is to obtain more compact representations of the data with limited
loss of information. Let $A$ be $m\times n$ matrix, then the low
rank approximation (rank $k$) of $A$ is given by 
\[
A_{m\times n}\approx B_{m\times k}C_{k\times n}.
\]
The low rank approximation of the matrix can be stored and manipulated
more economically than the matrix itself. One can see from the above
approximation that only $k(m+n)$ entries have to be stored instead
of $mn$ entries of the original matrix $A.$

The low rank approximation of a matrix appears in many applications.
The list of applications includes image processing \cite{FriedMehr,website},
data mining \cite{Elden,skillicorn}, noise reduction, seismic inversion,
latent semantic indexing \cite{website2}, principal component analysis
(PCA) \cite{pcabook,eldenpark}, machine-learning \cite{leekim,murphy,jye},
regularization for ill-posed problems, statistical data analysis applications,
DNA microarray data, web search model and so on. The low rank approximation
of matrices also plays a very important role in tensor decompositions
\cite{espigkumarjan,larsgras,boris1,boris2,boris3,boris4,lathauver,oseledets}.

Because of the interplay of rank and error there are basically two
types of problems related to the low rank approximation of a matrix;
fixed-precision approximation problem and fixed-rank approximation
problem (we do not use this nomenclature in this article). In the
fixed-precision approximation problem, for a given matrix $A$ and
a given tolerance $\epsilon,$ one wants to find a matrix $B$ with
rank $k=k(\epsilon)$ such that $\left\Vert A-B\right\Vert \leq\epsilon$
in an appropriate matrix norm. On the contrary, in the fixed-rank
approximation problem, one looks for a matrix $B$ with fixed rank
$k$ and an error $\left\Vert A-B\right\Vert $ as small as possible.

The low rank approximation problem is well studied in the numerical
linear algebra community. There are very classical matrix decompositions
which gives low rank approximation. Singular value decomposition (SVD)
is the best known. It has wide applications in many areas. It provides
the true rank and gives the best low rank approximation of a matrix
\cite{datta,golubvan,Sundarapan,trefthenbau}. QR decomposition with
column pivoting \cite{Stewartbook}, rank revealing QR factorization
\cite{tfchan,golubvan,Golubbusinger,golubstewartklema,guein} and
interpolative decomposition \cite{aricemgilaka,chengmartrokh,liberwoolfemartirokhlintygert,lucasmark}
are other useful techniques. These techniques require $O(mnk)$ arithmetic
operations to get a rank $k$ approximation by at least $k$ passes
(the number of times that the entire data is read) through the input
matrix. It is not easy to access the data in many applications with
very large data. So these methods become unsuitable for large scale
data matrices.

Alternatives for these classical algorithms are randomized algorithms
for low rank approximations \cite{deshpandevempala,Halkotropp,kannanvempala,MATINROHLINTYGERT,sarlos,sapp}.
The complexity of these algorithms is at most sublinear in the size
$m\times n$ and they only require one or two passes of the input
matrix. The main idea of these randomized algorithms is to compute
an approximate basis for the range space of the matrix $A$ using
a random selection of columns/rows of $A$ and project $A$ onto the
subspace spanned by this basis. We sketch it here:

Let $k$ be the target rank (the aim is to obtain a rank $k$ approximation).
Choose a number of samples larger than $k,$ i.e $s=k+p.$ The randomized
low rank approximation constructs the approximation in the following
way.

Step 1: Form lower dimensional matrix $X$ by the $s$ selected row
and/or columns.

Step 2: Compute an approximate orthonormal basis $Q=[q_{1},q_{2},...,q_{k}]$
for the range of $X$.

Step 3: Construct the low rank approximation $\tilde{A}$ by projecting
$A$ onto the space spanned by the basis $Q:\,\,\,\,\,\,\tilde{A}=QQ^{T}A.$

In step 1, the columns/rows can be chosen in different ways: by subsampling
of the input matrix or by using random projections. The matrix $X$
formed by these columns is expected to be very close to $A$ in a
sense that the basis of the range of $X$ covers the range of $A$
well. The orthonormal basis consisting of $k$ linearly independent
vectors can be obtained using exact methods since the size of $X$
is very small. These techniques are relatively insensitive to the
quality of randomness and produce high accurate results. The probability
of failure is negligible. Using the orthonormal basis $Q$ one can
approximate the standard factorizations like SVD, QR etc \cite{Halkotropp}.

There are other approximation techniques available in the literature
like cross/skeleton decompositions \cite{Bebendorf,Gantmacher,Goritryzamar,Trytyshivkov}.
Their complexity is of order $O(k^{2}(m+n))$ and they use only $k(m+n)$
entries from the original matrix to construct a rank $k$ approximation
of the matrix. These methods are also very useful in data sparse representation
of the higher order tensors. The algorithms that construct different
data tensor formats use low rank approximations of matrices at different
levels of their construction. These are obtained by the cross/skeleton
approximations inexpensively (linear in $m$ and $n$ ) which also
gives the data sparse representation with linear complexity \cite{espigkumarjan,boris1,oseledets}.

The main motivation of this paper is to give a brief description of
the techniques which are available in the literature. The paper gives
an overview of the existing classical deterministic algorithms, randomized
algorithms and finally cross/skeleton approximation techniques which
have great advantage of being able to handle really large data appearing
in applications.

In section 2 the classical algorithms like singular value decomposition,
pivoted QR factorization, rank revealing QR factorization (RRQR) are
described briefly with relevant references. More emphasize is given
to the subset selection problem and interpolative decomposition (these
play a big role in skeleton/cross approximation or $CUR$ decomposition
which will be discussed in section 3). Various randomized algorithms
are also described. In section 3 various versions of so called cross/skeleton
approximation techniques are described. The algorithms are given in
detail and the computational complexity of them are derived (linear
in $n$).

For simplicity of the presentation we consider only the matrix of
real numbers. Frobenius norm of a $m\times n$ matrix $A=(a_{ij})$
is defined as the square root of the sum of the absolute squares of
its elements i.e

\noindent 
\[
\left\Vert A\right\Vert _{F}=\sqrt{\sum_{i=1}^{m}\sum_{j=1}^{n}\left|a_{ij}\right|^{2}}.
\]

The spectral norm of the matrix $A$ is defined as largest singular
value of $A$. i.e 
\[
\left\Vert A\right\Vert _{2}=\sigma_{max}.
\]
 Here $\sigma_{max}$ is the largest singular value of the matrix
$A.$

\section{Classical techniques and randomized algorithms}

\subsection{Singular value decomposition}

This is a powerful technique in linear algebra. It gives the best
low rank approximation of a matrix. As mentioned in \cite{Stewartbook},
\textit{the singular value decomposition is the creme de la creme
of rank-reducing decompositions\textemdash the decomposition that
all others try to beat}.

Singular Value Decomposition factorizes $A\in\mathbb{R}^{m\times n}$
(where $m>n$), into the matrices $U$, $S$ and $V^{T}$, where $V^{T}$
is the transpose of a matrix $V$. The SVD factorization does not
require square matrices, therefore $m$, the number of rows, does
not have to equal $n$, the number of columns.

\noindent \emph{\noun{ 
\begin{equation}
A_{m\times n}=U_{m\times m}S_{m\times n}V_{n\times n}^{T}.\label{eq:}
\end{equation}
}}In this U and V are orthogonal matrices, therefore, all of the columns
of $U$ and $V$ are orthogonal to one another. The matrix $S$ is
a $m\times n$ rectangular diagonal matrix whose entries are in descending
order, $\sigma_{1}\geq\sigma_{2}\geq\cdots\geq\sigma_{n}\geq0,$ along
the main diagonal. 
\[
S=\left[\begin{array}{cccc}
\sigma_{1} & 0 & \cdots & 0\\
0 & \sigma_{2} & \cdots & 0\\
\vdots & \vdots & \ddots & 0\\
0 & 0 & \cdots & \sigma_{n}\\
0 & 0 & \cdots & 0\\
0 & 0 & \cdots & 0
\end{array}\right].
\]

\noindent \textbf{\textit{Note}}\textbf{:} If $A$ is a complex matrix,
then the singular value decomposition of a matrix $A$ is

\noindent 
\[
A_{m\times n}=U_{m\times m}S_{m\times n}V_{n\times n}^{*}
\]

\noindent where $V^{*}$ is the conjugate transpose of $V$ and $U$,
$V$ are unitary matrices.

\subsubsection*{Thin SVD }

Since $m>n,$ one can represent the SVD of $A$ as 
\[
A_{m\times n}=U_{m\times n}S_{n\times n}V_{n\times n}^{T}.
\]
Here $S_{n\times n}=diag(\sigma_{1},\sigma_{2},....,\sigma_{n}).$
This representation is called thin SVD of $A.$

\subsubsection*{Low-Rank Approximation}

\noindent Singular value decomposition gives the rank of a matrix.
The number of nonzero singular values of $A$ is the rank of the matrix
$A.$ Let the rank of $A$ be $r=min(m,n)$, then its SVD reads
\begin{eqnarray*}
A=U_{m\times r}S_{r\times r}V_{r\times n}^{T},\,\,\,\,\,\,\,\,\,\,\,\,\,\,\,\,\,\,\,\,\,\,\,\,\\
A=[u_{1}u_{2}....u_{r}]\left[\begin{array}{cccc}
\sigma_{1} &  &  & 0\\
 & \sigma_{2}\\
 &  & \ddots\\
0 &  &  & \sigma_{r}
\end{array}\right]\left[\begin{array}{c}
v_{1}^{T}\\
v_{2}^{T}\\
\vdots\\
v_{r}^{T}
\end{array}\right],\\
A=\sum_{i=1}^{r}\sigma_{i}u_{i}v_{i}^{T}\,\,\,\,\,\,\,\,\,\,\,\,\,\,\,\,\,\,\,\,\,\,\,\,\,\,\,\,\,\,\,\,\,
\end{eqnarray*}

\noindent where $u_{1},u_{2},..,u_{r}$ are columns of $U_{m\times r}$
and $v_{1},v_{2},...,v_{r}$ are columns of $V_{n\times r}.$ One
can see that the matrix $A$ is represented by the sum of outer products
of vectors. The matrix approximation by a low rank matrix is possible
using SVD. 

\noindent The rank $k$ approximation (also called as truncated or
partial SVD) of $A,$ $A_{k}$ where $k<r$, is given by zeroing out
the $r-k$ trailing singular values of $A$, that is 
\begin{align*}
A_{k} & =U_{m\times k}(S_{k})_{k\times k}V_{k\times n}^{T}=\sum_{i=1}^{k}\sigma_{i}u_{i}v_{i}^{T}.
\end{align*}

\noindent Here $S_{k}=diag(\sigma_{1},\sigma_{2},....,\sigma_{k}),\,U_{m\times k}=[u_{1},u_{2},..,u_{k}]\,\,\textrm{and}\,\,V_{k\times n}^{T}=\left[\begin{array}{c}
v_{1}^{T}\\
v_{2}^{T}\\
\vdots\\
v_{k}^{T}
\end{array}\right]$. Then one can see that 
\begin{eqnarray*}
A_{k}=U_{m\times k}U_{k\times m}^{T}A=\left(\sum_{i=1}^{k}u_{i}u_{i}^{T}\right)A\,\,\,and\,\,\,A_{k}=AV_{n\times k}V_{k\times n}^{T}=A\left(\sum_{i=1}^{k}v_{i}v_{i}^{T}\right),
\end{eqnarray*}
 i.e $A_{k}$ is the projection of the $A$ onto the space spanned
by the top $k$ singular vectors of $A.$ The following theorem states
that the above approximation is the best rank $k$ approximation in
both Frobenius and spectral norm.

\subsubsection*{Theorem: Eckart-Young theorem}

Let $A_{k}$ be the rank-$k$ approximation of $A$ achieved by SVD-truncation
as above. Then $A_{k}$ is the closest rank-$k$ matrix to $A$, i.e.
\begin{eqnarray}
\left\Vert A-A_{k}\right\Vert _{F}\leq\left\Vert A-B\right\Vert _{F}\label{eq:-1}
\end{eqnarray}
where $B$'s are rank-$k$ matrices.

\noindent The minimal error is given by the Euclidean norm of the
singular values that have been zeroed out in the process 
\[
\left\Vert A-A_{k}\right\Vert _{F}=\sqrt{\sigma_{k+1}^{2}+\cdots+\sigma_{r}^{2}}
\]
where $\left\Vert .\right\Vert _{F}$ is Frobenius norm.\medskip{}

\noindent \textbf{Remark}: SVD also gives the best low rank approximation
in spectral norm: 
\begin{eqnarray*}
\left\Vert A-A_{k}\right\Vert _{2}=\min_{rank(B)=k}\left\Vert A-B\right\Vert _{2}=\sigma_{k+1}.
\end{eqnarray*}

\subsubsection*{Algorithms and computational complexity}

\noindent The SVD of a matrix $A$ is typically computed numerically
by a two-step procedure. In the first step, the matrix is reduced
to a bidiagonal matrix. This takes $O(mn^{2})$ floating-point operations
(flops), and the second step is to compute the SVD of the bidiagonal
matrix. The second step takes $O(n)$ iterations, each costing $O(n)$
flops. Therefore the overall cost is $O(mn^{2})$. If $A$ is a square
matrix, then SVD algorithm requires $O(n^{3})$ flops \cite{golubvan,trefthenbau}.

Alternatively we can obtain the rank $k$ approximation directly by
obtaining partial SVD. The partial SVD can be obtained by computing
partial QR factorization and post process the factors \cite{Halkotropp}.
This technique requires only $O(kmn)$ flops. Krylov subspace methods
like, Lanczos methods for certain large sparse symmetric matrices
and Arnoldi (unsymmetric Lanczos methods) for unsymmetric matrices
can be used to compute SVD \cite{golubvan}. The straight forward
algorithm of Lanczos bidiagonalization has the problem of loss of
orthogonality between the computed Lanczos vectors. Lanczos with complete
reorthogonalization (or performing only local orthogonalization at
every Lanczos steps), block Lanczos algorithms are practical Lanczos
procedures \cite{golubvan}. Details of efficient algorithms for large
sparse matrices can be found in \cite{lanc5} , chapter 4 of \cite{lanc6}
and \cite{lanc3,lanc8,lanc4}. The algorithms are available in the
packages SVDPACK \cite{lanc1,lanc7}, PROPACK \cite{lanc2}. 

\medskip{}

As a low rank approximation method the singular value decomposition
has few drawbacks. It is expensive to compute if the dimension of
the matrix is very large. In many applications it is sufficient to
have orthonormal bases for the fundamental subspaces, something which
the singular value decomposition provides. In other applications,
however, it is desirable to have natural bases that consist of the
rows or columns of the matrix. Here we describe such matrix decompositions.

\subsection{Pivoted QR decomposition}

Let $A_{m\times n}$ be a rank deficient matrix ($m>n$) with rank
$\gamma$. A pivoted QR decomposition with column pivoting has the
form 
\begin{eqnarray*}
 &  & AP=QR
\end{eqnarray*}
where $P$ is a permutation matrix, $Q$ is orthonormal and $R$ is
upper triangular matrix. In exact arithmetic, 
\begin{eqnarray*}
AP=Q\left[\begin{array}{cc}
R_{11}^{(\gamma)} & R_{12}^{(\gamma)}\\
0 & 0
\end{array}\right]
\end{eqnarray*}
where $R_{11}^{(\gamma)}$ is $\gamma\times\gamma$ upper triangular
matrix with rank $\gamma,$ $Q\in\mathbb{R}^{m\times n}$ and $P\in\mathbb{R}^{n\times n}.$
In floating point arithmetic one may obtain 
\begin{eqnarray}
AP=Q\left[\begin{array}{cc}
R_{11}^{(\gamma)} & R_{12}^{(\gamma)}\\
0 & R_{22}^{(\gamma)}
\end{array}\right]\label{eq:3}
\end{eqnarray}
such that $\left\Vert R_{22}^{(\gamma)}\right\Vert $ is small.

A rank $k$ approximation to any matrix $A$ can be obtained by partitioning
the decomposition $AP=QR.$ Let $B=AP$ and write 
\begin{eqnarray}
B=[B_{1}^{(k)}B_{2}^{(k)}]=[Q_{1}^{(k)}Q_{2}^{(k)}]\left[\begin{array}{cc}
R_{11}^{(k)} & R_{12}^{(k)}\\
0 & R_{22}^{(k)}
\end{array}\right],\label{eq:4}
\end{eqnarray}
where $B_{1}^{(k)}$ has $k$ columns. Then our rank $k$ approximation
is 
\begin{eqnarray*}
\hat{B}^{(k)}=Q_{1}^{(k)}[R_{11}^{(k)}R_{12}^{(k)}].
\end{eqnarray*}
Therefore 
\begin{eqnarray*}
B-\hat{B}^{(k)}=Q_{2}^{(k)}[0\,\,R_{22}^{(k)}].
\end{eqnarray*}
The approximation $\hat{B}^{(k)}$ reproduces the first $k$ columns
of $B$ exactly. Since $Q_{2}^{(k)}$ is orthogonal, the error in
$\hat{B}^{(k)}$ as an approximation to $B$ is 
\begin{eqnarray}
\left\Vert B-\hat{B}^{(k)}\right\Vert =\left\Vert R_{22}^{(k)}\right\Vert .\label{eq:5}
\end{eqnarray}

This is called truncated pivoted QR decomposition to $A.$ The permutation
matrix is determined by column pivoting such that $R_{11}^{(k)}$
is well conditioned and $R_{22}^{(k)}$ is negligible (the larger
entries of $R$ are moved to the upper left corner and the smallest
entries are isolated in the bottom submatrix). This decomposition
is computed by a variation of orthogonal triangularization by Householder
transformations \cite{golubvan,Golubbusinger,Stewartbook}. The algorithm
described in \cite{golubvan,Stewartbook} requires $O(kmn)$ flops.
This algorithm is effective in producing a triangular factor $R$
with small $\left\Vert R_{22}^{(k)}\right\Vert ,$ very little is
known in theory about its behavior and it can fail on some matrices
(look at example 1 of \cite{tfchan} and also in \cite{faddeevkub}
).

Similar decompositions like pivoted Cholesky decompositions, pivoted
QLP decomposition and UTV decompositions can be found in \cite{Stewartbook}.

\subsubsection*{Rank revealing QR factorization}

As we have seen above, the column pivoting QR decomposition is a cheaper
alternative to SVD. This factorization works in many cases but may
also fail sometimes \cite{tfchan,faddeevkub}. The most promising
alternative to SVD is the so-called rank revealing QR factorization.
\medskip{}
\\
 \textbf{Definition (RRQR):} Given a matrix $A_{m\times n}(m\geq n)$
and an integer $k$($k\leq n),$ assume partial QR factorizations
of the form 
\begin{eqnarray}
AP=QR=Q\left[\begin{array}{cc}
R_{11} & R_{12}\\
0 & R_{22}
\end{array}\right],\label{eq:6}
\end{eqnarray}
where $Q\in\mathbb{R}^{m\times n}$ is an orthonormal matrix, $R\in\mathbb{R}^{n\times n}$
is block upper triangular, $R_{11}\in\mathbb{R}^{k\times k}$, $R_{12}\in\mathbb{R}^{k\times n-k},$~$R_{22}\in\mathbb{R}^{n-k\times n-k}$
and $P\in\mathbb{R}^{n\times n}$ is a permutation matrix. The above
factorization is call RRQR factorization if it satisfies 
\begin{eqnarray}
\frac{\sigma_{k}(A)}{p(k,n)}\leq\sigma_{min}(R_{11})\leq\sigma_{k}(A)\label{7(a)}\\
\sigma_{k+1}(A)\leq\sigma_{max}(R_{22})\leq p(k,n)\sigma_{k+1}(A)\label{7(b)}
\end{eqnarray}
where $p(k,n)$ is a lower degree polynomial in $k$ and $n.$

In the above definition $\sigma_{min}$ is the minimum singular value
and $\sigma_{max}$ is the maximum singular value. RRQR was defined
by Chan in \cite{tfchan}(similar ideas were proposed independently
in \cite{Foster}). A constructive proof of the existence of a RRQR
factorization of an arbitrary matrix $A_{m\times n}$ with numerical
rank $r$ is given in \cite{honpan}. Much research on RRQR factorizations
has yielded improved results for $p(k,n).$ There are several algorithms
to compute the RRQR factorization \cite{tfchan,chandraispin,golubvan,guein}.
The computational complexity of these algorithms are slightly larger
than the standard QR decomposition algorithm. The values of $p(k,n)$
and the complexities of different algorithms were tabulated in \cite{boumahodri}.

Different applications of RRQR like subset selection problems, total
least-squares problems including low rank approximation have been
discussed in \cite{chanhansen}. The low rank approximation of the
matrix $A$ can be obtained by neglecting the submatrix $R_{22}$
in RRQQ factorization of $A.$ It has been shown that matrix approximations
derived from RRQR factorizations are almost as good as those derived
from truncated SVD approximations.

The singular value and pivoted QR decompositions are not good for
large and sparse matrices. The problem is that the conventional algorithms
for computing these decomposition proceed by transformations that
quickly destroy the sparsity of matrix $A$. Different algorithms
for the efficient computation of truncated pivoted QR approximations
to a sparse matrix without loosing the sparsity of the matrix $A$
are proposed in \cite{berrypulotovastewart,stewartpaper1}. Some more
references on structure preserving RRQR factorization algorithms are
given in \cite{chanhansen}.

\subsection{Interpolative decomposition}

Interpolative decompositions (ID's) (also called as $CX$ decomposition)
are closely related to pivoted QR factorizations and are useful for
representing low rank matrices in terms of linear combinations of
their columns \cite{chengmartrokh,liberwoolfemartirokhlintygert,MATINROHLINTYGERT}.
Interpolative decomposition of a matrix completely rely on the column
subset selection. Before defining the interpolative decomposition,
a brief description is given below on the subset selection problem.

\subsubsection*{Subset selection problem}

Subset selection is a method for selecting a subset of columns from
a real matrix, so that the subset represents the entire matrix well
and is far from being rank deficient. Given a $m\times n$ matrix
$A$ and an integer $k,$ subset selection attempts to find the $k$
most linearly independent columns that best represents the information
in the matrix.

The mathematical formulation of the subset selection problems is:
Determine a permutation matrix $P$ such that 
\begin{equation}
AP=(A_{1}\,A_{2}),\,\textrm{where}\label{eq:9}
\end{equation}
~~~~1. $A_{1}$ is $m\times k$ matrix containing $k$ linearly
independent columns such that smallest singular value is as large
as possible. That is for some $\gamma$ 
\begin{eqnarray}
\frac{\sigma_{k}(A)}{\gamma}\leq\sigma_{k}(A_{1})\leq\sigma_{k}(A).\label{eq:-2}
\end{eqnarray}

2. the $n-k$ columns of $A_{2}$ (redundant columns) are well represented
by $k$ columns of $A_{1}.$ That is 
\begin{eqnarray}
\begin{array}{c}
min\\
Z\in\mathbb{R}^{k\times n-k}
\end{array}\left\Vert A_{1}Z-A_{2}\right\Vert _{2}\,\,\,\,\,\textrm{is\,\,\ small. }\nonumber \\
\textrm{i.e for\,\,\ some}\,\,\gamma,\,\,\,\,\,\,\,\sigma_{k+1}(A)\leq\begin{array}{c}
min\\
Z\in\mathbb{R}^{k\times n-k}
\end{array}\left\Vert A_{1}Z-A_{2}\right\Vert _{2}\,\leq\gamma\,\,\sigma_{k+1}(A).\label{eq:-3}
\end{eqnarray}
\textbf{Remark:} \textit{$Z$ is a matrix responsible for representing
the columns of $A_{2}$ in terms of the columns of $A_{1}.$ \smallskip{}
 }

More detailed information and an equivalent definition to the subset
selection problem is given in section 2.4. \textit{ }The subset selection
using singular value decomposition has been addressed in \cite{golubvan,golubstewartklema}.
Many subset selection algorithms use a QR decomposition (which was
discussed in last subsection) to find the most representative columns
\cite{brodbrownpenne}. There are several randomized algorithms for
this problem \cite{boumahodri,boumahodri2,civrilismail}. The strong
RRQR algorithm by Gu and Eisenstat \cite{guein} gives the best deterministic
approximation to the two conditions (10) and (11) of the subset selection
problem. The details are given below.

As described in the last subsection the RRQR factorization of $A_{m\times n}$
is represented by 
\begin{eqnarray*}
AP=QR=Q\left[\begin{array}{cc}
R_{11} & R_{12}\\
0 & R_{22}
\end{array}\right].
\end{eqnarray*}

This gives the permutation matrix $P$ such that $AP=(A_{1},A_{2})$
where $A_{1}$(the matrix with most important $k$ columns of $A$)
and $A_{2}$ (with the redundant columns) are given by 
\begin{eqnarray*}
A_{1}=Q\left[\begin{array}{c}
R_{11}\\
0
\end{array}\right]\,\,\,\textrm{and}\,\,\,A_{2}=Q\left[\begin{array}{c}
R_{12}\\
R_{22}
\end{array}\right]
\end{eqnarray*}
with 
\begin{eqnarray*}
\sigma_{i}(A_{1})\geq\frac{\sigma_{i}(A)}{\sqrt{1+f^{2}k(n-k)}},1\leq i\leq k\\
\begin{array}{c}
\textrm{min}\\
Z
\end{array}\,\left\Vert A_{1}Z-A_{2}\right\Vert _{2}\leq\sigma_{k+1}(A)\sqrt{1+f^{2}k(n-k)}.
\end{eqnarray*}

In the above inequalities $f\geq1$ is a tolerance supplied by some
user. The Gu and Eisenstat algorithm also guarantees that $\left|(R_{11}^{-1}R_{12})_{ij}\right|\leq1$
for $1\leq i\leq k,1\leq j\leq n-k.$ One can extend this algorithm
for wide and fat matrices where $m<n$ and $k=m.$ The computational
complexity of this algorithm is $O(mn^{2}).$ \medskip{}
 \\
 \textbf{Remark: }\textit{So from the strong RRQR algorithm we can
see the following.}

\textit{As described in subsection 2.2 the truncated RRQR of $A$
is $AP\simeq Q[R_{11}\,\,R_{12}].$ Now we can write it as $AP\simeq QR_{11}[I_{k\times k}\,\,R_{11}^{-1}R_{12}],$
where $QR_{11}$ is matrix which contains $k$ linearly independent
columns. From theorem 3.2 of \cite{guein}, one can see that $\left\Vert [I\,\,\,R_{11}^{-1}R_{12}]\right\Vert _{2}\leq\sqrt{1+f^{2}k(n-k)}.$
Therefore $QR_{11}[I\,\,\,R_{11}^{-1}R_{12}]P^{T}$ is an approximation
to the matrix $A.$ As we have seen in section 2.2 the error in the
approximation is 
\begin{eqnarray*}
\left\Vert A-QR_{11}[I\,\,R_{11}^{-1}R_{12}]P^{T}\right\Vert _{2}=\left\Vert R_{22}\right\Vert _{2}\leq\sigma_{k+1}(A)\,\sqrt{1+f^{2}k(n-k)}.
\end{eqnarray*}
 \medskip{}
}

This subset selection problem is also studied widely in randomized
setting. We postpone the discussion of these techniques to subsection
2.4.\medskip{}
\\
\textbf{Remark:} \textit{In most of the deterministic algorithms for
subset selection problem, the error estimates are given in spectral
norm. The error estimates in both spectral and Frobenius norm are
presented for several randomized algorithms in the literature which
is the subject of next subsection.} \medskip{}
\\
\textbf{Definition (ID): }Let $A_{m\times n}$ be a matrix of rank
$k.$ There exists an $m\times k$ matrix $B$ whose columns constitute
a subset of the columns of $A,$ and $k\times n$ matrix $P,$ such
that

1. some subset of the columns of $P$ makes up $k\times k$ identity
matrix,

2. $P$ is not too large (no entry of $P$ has an absolute value greater
than $1$), and

3. $A_{m\times n}=B_{m\times k}P_{k\times n}.$

Moreover, the decomposition provides an approximation 
\[
A_{m\times n}\approx B_{m\times k}P_{k\times n}
\]
when the exact rank of $A$ is greater than $k,$ but the $(k+1)st$
greatest singular value of $A$ is small. The approximation quality
of the Interpolative decomposition is described in the following Lemma
\cite{liberwoolfemartirokhlintygert,MATINROHLINTYGERT}. One can also
look at \cite{GoreTryty,gorotyrzama1} for similar results. \smallskip{}
 \\
 \textbf{Lemma:} Suppose that $m$ and $n$ are positive integers,
and $A$ is $m\times n$ matrix. Then for any positive integer $k$
with $k\leq m$ and $k\leq n,$ there exist a $k\times n$ matrix
$P,$ and a $m\times k$ matrix $B$ whose columns constitute a subset
of the columns of $A,$ such that

1. some subset of the columns of $P$ makes up $k\times k$ identity
matrix,

2. no entry of $P$ has an absolute value greater than $1,$

3. $\left\Vert P_{k\times n}\right\Vert _{2}\leq\sqrt{k(n-k)+1}$,

4. the least (that is, $k$ the greatest) singular value of $P$ is
at least 1,

5. $A_{m\times n}=B_{m\times k}P_{k\times n}$ when $k=m$ and $k=n,$
and

6. when $k<m$ and $k<n,$ 
\begin{eqnarray*}
\left\Vert A_{m\times n}-B_{m\times k}P_{k\times n}\right\Vert _{2}\leq\sqrt{k(n-k)+1}\,\,\sigma_{k+1},
\end{eqnarray*}
~~~~~~~where $\sigma_{k+1}$ is the $(k+1)st$ greatest singular
value of $A.$

The algorithms to compute the ID are computationally expensive. The
algorithms described in \cite{guein} can be used to compute the Interpolative
decomposition. In \cite{MATINROHLINTYGERT} a randomized algorithm
has been proposed. The authors have constructed the interpolative
decomposition under weaker conditions than those in above Lemma. The
computational complexity of this algorithm is $O(kmnlog(n)).$ This
decomposition have also been studied in \cite{aricemgilaka,lucasmark}.
The details of a software package of ID algorithms can be found in
\cite{martintygertsoftware} and the applications of ID in different
applications can be found in \cite{comptonosher,panweipeng}.

\subsubsection*{CUR decomposition}

Interpolative decomposition can also be used to obtain the independent
rows of a matrix. So two ID's can be combined to construct the matrix
using the subcollection of its columns and rows which is called CUR
decomposition.

A Given matrix $A_{m\times n}$ is decomposed into a product of three
matrices $C,U$ and $R,$ where $C$ consists of small number of actual
columns of $A,$ $R$ consists of a small number of actual rows of
$A$ and $U$ is a small carefully constructed matrix that guarantees
that the product $CUR$ is close to $A.$

This decomposition is also known as skeleton decomposition \cite{Gantmacher}.
Let $A_{m\times n}$ be a matrix with rank $r.$ One can reconstruct
the matrix $A$ by choosing $C_{m\times r}$ with $r$ columns and
$R_{r\times n}$ with $r$ rows of the matrix such that the intersection
matrix $W_{r\times r}$ is nonsingular, the corresponding $CUR$ decomposition
is 
\[
A=CUR,\,\,\,\textrm{with}\,\,\,U=W^{-1}.
\]
This decomposition provides a low rank approximation $A_{k}\simeq CUR$
when the number of selected rows and columns is less than $r.$ The
proper choice of these columns and rows and the matrix $U$ was the
subject of study in the literature.

In \cite{berrypulotovastewart,stewartpaper1} a Gram-Schmidt based
pivoted QR algorithm is proposed. The matrices $C$ and $R$ were
obtained by applying the algorithm to the matrices $A,$ $A^{T}$
respectively. The matrix $U$ has been obtained such that (see section
4 of \cite{berrypulotovastewart}) 
\[
\left\Vert A-CUR\right\Vert _{F}^{2}=\textrm{min}.
\]

In \cite{GoreTryty,Goritryzamar} a $CUR$ decomposition has been
developed which is called pseudoskeleton approximation. The $k$ columns
in $C$ and $k$ rows in $R$ were chosen such that their intersection
$W_{k\times k}$ has maximum volume (the maximum determinant among
all $k\times k$ submatrices of $A$). The complete description of
this class of algorithms is given in Section 3. Randomized algorithms
for $CUR$ decomposition are discussed in the following subsection.

\subsection{Randomized algorithms}

There are several randomized algorithms to obtain the low rank approximation
of a matrix in the literature. As described in section 1 (step 1),
a lower dimensional matrix $X$ can be constructed by selecting $s$
columns/rows. These columns/rows can be chosen in different ways by
subsampling of the given matrix or by random projection. First we
look at subsampling based randomized algorithms for the subset selection
problem and for the $CUR$ decomposition and we discuss the random
projection method. We also discuss randomized SVD. \smallskip{}

\subsubsection*{(I). Sampling based methods\smallskip{}
}

\subsubsection*{Column subset selection problem (CSSP)}

We have seen the mathematical formulation of the subset selection
problem in section 2.3. There is an equivalent definition to the subset
selection problem in the literature. Here we define the column subset
selection problem.

\textbf{$\!\!\!\!\!\!\!\!$Definition (CSSP):} Given a matrix $A_{m\times n}$
and a positive integer $k$ as the number of columns of $A$ forming
a matrix $C\in\mathbb{R}^{m\times k}$ such that the residual $\left\Vert A-P_{C}A\right\Vert _{\xi}$
is minimized over all possible $\left(\begin{array}{c}
n\\
k
\end{array}\right)$ choices for the matrix $C.$ Here $P_{C}=CC^{\dagger}$ ($C^{\dagger}$
is Moore-Penrose pseudoinverse of $C$) denotes the projection onto
the $k$ dimensional space spanned by the columns of $C$ and $\xi=2\,\,or\,\,F$
denotes the spectral norm or Frobenius norm.

This seems to be a very hard optimization problem, finding $k$ columns
out of $n$ columns such that $\left\Vert A-P_{C}A\right\Vert _{\xi}$
is minimum. It requires $O(n^{k})$ time and thus to find the optimal
solution we require $O(n^{k}mnk).$ So obtaining the approximation
is prohibitively slow if the data size is large. The NP-hardness of
the CSSP (assuming $k$ is a function of $n$) is an open problem
\cite{boumahodri2}. So research is focused on computing approximation
solutions to CSSP.

Let $A_{k}$ be the best low rank $k$ approximation. Therefore $\left\Vert A-A_{k}\right\Vert _{\xi}$
provides a lower bound for $\left\Vert A-P_{C}A\right\Vert _{\xi}$
for $\xi=F,2$ and for any choice of $C.$ So most of the algorithms
have been proposed in the literature to select $k$ columns of $A$
such that the matrix $C$ satisfies 
\begin{eqnarray*}
\left\Vert A-A_{k}\right\Vert _{\xi}\leq\left\Vert A-P_{C}A\right\Vert _{\xi}\leq p(k,n)\left\Vert A-A_{k}\right\Vert _{\xi}
\end{eqnarray*}
for some function $p(k,n).$

As we have seen in the previous section, the strong RRQR algorithm
(deterministic algorithm) gives spectral norm bounds. From the definition
of RRQR there exists a permutation matrix $\Pi\in\mathbb{R}^{n\times n}$
(look at equation (6), please note that symbol for permutation matrix
is changed here). Let $\Pi_{k}$ denote the first $k$ columns of
this permutation matrix $\Pi.$ If $C=A\Pi_{k}$ is $m\times k$ matrix
consisting of $k$ columns of $A$ ($C$ corresponds to $Q\left[\begin{array}{c}
R_{11}\\
0
\end{array}\right]$ in definition of $RRQR$), then from the equations (5) and (8) one
can see (proof is very simple and similar to the proof of Lemma 7.1
in \cite{brodbrownpenne} 
\begin{eqnarray*}
\left\Vert A-P_{C}A\right\Vert _{2}=\left\Vert R_{22}\right\Vert _{2}\leq p(k,n)\sigma_{k+1}(A).
\end{eqnarray*}
That is, any algorithm that constructs an RRQR factorization of the
matrix A with provable guarantees also provides provable guarantees
for the CSSP \cite{boumahodri}.

Several randomized algorithms have been proposed to this problem.
In these methods, few columns (more than the target rank $k$) of
$A$ are selected randomly according to a probability distribution
which was obtained during the preprocessing of the matrix and then
the low rank approximation is obtained using classical techniques
from linear algebra. One such type of method, a fast Monte-Carlo algorithm
for finding a low rank approximation, has been proposed in \cite{frikannanvempala,frikannanvempala2}.
This algorithm gives an approximation very close to SVD by sampling
the columns and rows of the matrix $A$ with only two passes through
the data. It is based on selecting a small subset of important columns
of $A,$ forming a matrix $C$ such that the projection of $A$ on
the subspace spanned by the columns of $C$ is as close to $A$ as
possible. A brief description of the algorithm is given below.

A set of $s$ columns $(s>k,$ where $k$ is the target rank) were
chosen randomly, each according to a probability distribution proportional
to their magnitudes (squared $l_{2}$ norms of the columns). Let $S$
be the matrix obtained by writing these $s$ columns as columns. An
orthogonal set of $k$ vectors in the span of these $s$ columns have
been obtained. These orthogonal vectors are the top $k$ left singular
vectors of the matrix $S$ (which were obtained from the SVD of a
$s\times s$ matrix formed by sampling the rows according to a probability
distribution). The rank $k$ approximation to $A$ is obtained by
projecting $A$ on the span of these orthogonal vectors.

The rank $k$ approximation $D^{*}$ of the matrix $A$ (within a
small additive error) may be computed such that 
\begin{equation}
\left\Vert A-D^{*}\right\Vert _{F}^{2}\leq\left\Vert A-A_{k}\right\Vert _{F}^{2}+\epsilon\left\Vert A\right\Vert _{F}^{2},\label{eq:-4}
\end{equation}
holds with probability at least $1-\delta.$ Here $\delta$ is the
failure probability, $\epsilon$ is an error parameter and the randomly
chosen columns $s=poly(k,1/\epsilon)$ (a polynomial in $k$ and $1/\epsilon$).
Here $A_{k}$ denotes the best rank $k$ approximation of $A.$ It
requires $O(ms^{2}+s^{3})$ complexity, where $s=O\,\left(max(k^{4}\epsilon^{-2},k^{2}\epsilon^{-4})\right)$.
The matrix $D^{*}$ can be explicitly constructed in $O(kmn)$ time.
The additive error $\epsilon\left\Vert A\right\Vert _{F}^{2}$ could
be arbitrarily large (for the matrices with sufficiently large $\left\Vert A\right\Vert _{F}^{2}$)
compared to error $\left\Vert A-A_{k}\right\Vert _{F}^{2}.$ This
kind of sampling method may not perform well in some cases \cite{deshpvempalawang}.

In \cite{drikannamahoney} a modified version of the algorithm proposed
in \cite{frikannanvempala,frikannanvempala2} has been discussed.
In \cite{deshpvempalawang} Deshpande et. al. generalized the work
in \cite{frikannanvempala,frikannanvempala2}. They have proved that
the additive error in (12) drops exponentially by adaptive sampling
and presented a multipass algorithm for low rank approximation. They
have shown that it is possible to get $(1+\epsilon)$ relative or
multiplicative approximation (look at (14)). They have generalized
the sampling approach using volume sampling (i.e picking $k$-subsets
of the columns of any given matrix with probabilities proportional
to the squared volumes of the simplicies defined by them) to get a
multiplicative approximation (look at (14)) instead of additive approximation
(look at (12)). They have proved the following existence result. There
exists (using volume sampling) exactly $k$ columns in any $m\times n$
matrix $A$ such that 
\begin{eqnarray}
\left\Vert A-D^{*}\right\Vert _{F}\leq\sqrt{(k+1)}\left\Vert A-A_{k}\right\Vert _{F},\label{eq:14}
\end{eqnarray}
 where $D^{*}$ (this may not coincide with the $D^{*}$ in (12))
is the projection onto the span of these $k$ columns. 

They also have proved (existence result) that there exists $k+k(k+1)/\epsilon$
rows whose span contains the rows of a rank-$k$ matrix $D^{*}$ such
that 
\begin{equation}
\left\Vert A-D^{*}\right\Vert _{F}^{2}\leq(1+\epsilon)\,\,\left\Vert A-A_{k}\right\Vert _{F}^{2}.\label{eq:13}
\end{equation}

In \cite{deshpandevempala}, Deshpande et. al. improved the existence
result in eq (14) and developed an efficient algorithm. They have
used an adaptive sampling method to approximate the volume sampling
method and developed an algorithm which finds $k$ columns of $A$
such that 
\begin{eqnarray}
\left\Vert A-D^{*}\right\Vert _{F}\leq\sqrt{(k+1)!}\left\Vert A-A_{k}\right\Vert _{F}.\label{eq:15}
\end{eqnarray}
The computational complexity of this algorithm is $O(mnk+kn).$ This
algorithm requires multipasses through the data and also maintains
the sparsity of $A$. 

In \cite{boumahodri,boumahodri2}, Boutsidis et. al. proposed a two
stage algorithm to select exactly $k$ columns from a matrix. In the
first stage (randomized stage), the algorithm randomly selects $O(k\,ln\,k)$
columns of $V_{k}^{T},$ i.e of the transpose of the $n\times k$
matrix consisting of the top $k$ right singular vectors of $A,$
according to a probability distribution that depends on information
in the top-$k$ right singular subspace of $A.$ Then in the second
stage (the deterministic stage), $k$ columns have been selected from
the set of columns of $V_{k}^{T}$ using deterministic column selection
procedure. The computational complexity of this algorithm is $O(min(mn^{2},m^{2}n)).$
It has been proved that the algorithm returns a $m\times k$ matrix
$C$ consisting of exactly $k$ columns of $A$ (rank of $A$ is $\rho$)
such that with probability at least $0.7:$ 
\begin{eqnarray*}
\left\Vert A-P_{C}A\right\Vert _{2}\leq O\left(k^{3/4}\,log^{1/2}\,(k)\,(\rho-k)^{1/4}\right)\,\,\left\Vert A-A_{k}\right\Vert _{2},\\
\left\Vert A-P_{C}A\right\Vert _{F}\leq O\left(k\sqrt{log\,k}\right)\,\,\left\Vert A-A_{k}\right\Vert _{F}.\,\,\,\,\,\,\,\,\,\,\,\,\,\,\,\,\,\,\,\,\,\,\,\,\,\,\,\,\,\,\,\,\,\,\,
\end{eqnarray*}

They have compared the approximation results with best existing results
for CSSP. They have shown that (from the above equations) the estimate
in spectral norm is better than the existing result ($\left\Vert A-P_{C}A\right\Vert _{2}\leq O(\sqrt{1+k(n-k)})\left\Vert A-A_{k}\right\Vert _{2}$~~\cite{guein})
by a factor of $n^{1/4}$ and worse than the best existence result
(13) by a factor $O(\sqrt{k\,log\,k})$ in Frobenius norm. 

In \cite{deshpanderadmacher} Deshpande and Rademacher have proposed
an efficient algorithm for volume sampling. They have selected $k$
columns such that 
\begin{eqnarray*}
\left\Vert A-P_{C}A\right\Vert _{2}\leq\sqrt{(k+1)(n-k)}\left\Vert A-A_{k}\right\Vert _{2}\\
\left\Vert A-P_{C}A\right\Vert _{F}\leq\sqrt{(k+1)}\left\Vert A-A_{k}\right\Vert _{F}\,\,\,\,\,\,\,\,\,\,\,
\end{eqnarray*}
with $O(kmn^{\omega}\,log\,n)$ arithmetic operations ($\omega$ is
the exponent of arithmetic complexity of matrix multiplication). This
improves the $O(k\sqrt{log\,k})-$ approximation of Boutsidis et.
al. in \cite{boumahodri} for the Frobenius norm case.

In the very recent articles by Boutsidis et. al \cite{boudriismail}
and Guruswami et.al. \cite{guruswamisinop}, these estimates have
been further improved. This problem has been also studied in \cite{haraimaung,civrilismail,drineasmahoneymuthu,farahatghodsi,friedzare,harpeled,mahoney,MAUNGSCH,pipezhouzhang,rudelsonver,wiitencandes,woolflibertyrokh}
and also in a PhD thesis by Civril \cite{civrilthesis}. In \cite{clarkwood,GHASAMI,ghashami2,liberty,woodruff},
a similar kind of work has been studied (streaming algorithms).

\subsubsection*{Randomised CUR }

As described in section 2.3, $CUR$ decomposition gives low rank approximation
explicitly expressed in terms of a small number of columns and rows
of the matrix $A.$ $CUR$ decomposition problem has been widely studied
in the literature. This problem has a close connection with the column
subset selection problem. One can obtain the CUR decomposition by
using column subset selection on $A$ and on $A^{T}$ to obtain the
matrices $C$ and $R$ respectively. But this will double the error
in the approximation. Most of the existing CUR algorithms uses column
subset selection procedure to choose the matrix $C.$

In \cite{drikannamahoney2}, Drineas et. al. have proposed a linear
time algorithm to approximate the CUR decomposition. $c$ columns
of $A$ and $r$ rows of $A$ are randomly chosen according to a probability
distribution to obtain the matrices $C_{m\times c}$ consisting of
chosen $c$ columns, $R_{r\times n}$ consisting of chosen $r$ rows.
A $c\times r$ matrix $U$ has been obtained using $C$ and $R.$
They have shown that for given $k$, by choosing $O(log(1/\delta)\epsilon^{-4})$
columns of $A$ to construct $C$ and $O(k\delta^{-2}\epsilon^{-2})$
rows of $A$ to construct $R$, the resulting CUR decomposition satisfies
the additive error bound with probability at least $1-\delta$ 
\begin{eqnarray*}
\left\Vert A-CUR\right\Vert _{2}\leq\left\Vert A-A_{k}\right\Vert _{2}+\epsilon\left\Vert A\right\Vert _{2} & .
\end{eqnarray*}
By choosing $O(klog(1/\delta)\epsilon^{-4})$ columns of $A$ to construct
$C$ and $O(k\delta^{-2}\epsilon^{-2})$ rows of $A$ to construct
$R$, the resulting CUR decomposition satisfies the additive error
bound with probability at least $1-\delta$

\begin{eqnarray*}
\left\Vert A-CUR\right\Vert _{F}\leq\left\Vert A-A_{k}\right\Vert _{F}+\epsilon\left\Vert A\right\Vert _{F} & .
\end{eqnarray*}
Here $\epsilon$ is the error parameter and $\delta$ is the failure
probability. The complexity of the algorithm is $O(mc^{2}+nr+c^{2}r+c^{3}),$
which is linear in $m$ and $n.$ This algorithm needs very large
number of rows and columns to get good accuracy.

In \cite{drimahomuthu}, Drineas et. al. developed an improved algorithm.
$c$ columns and $r$ rows were chosen randomly by subsampling to
construct the matrices $C$ and $R$ respectively and $U$ is the
weighted Moore-Penrose inverse of the intersection between the matrices
$C$ and $R.$ For given $k$, they have shown that there exists randomized
algorithms such that exactly $c=O(k^{2}log(1/\delta)\epsilon^{-2})$
columns of $A$ are chosen to construct $C,$ then exactly $r=O(c^{2}log(1/\delta)\epsilon^{-2})$
rows of $A$ are chosen to construct $R$, such that with probability
at least $1-\delta,$ 
\begin{eqnarray*}
\left\Vert A-CUR\right\Vert _{F}\leq(1+\epsilon)\left\Vert A-A_{k}\right\Vert _{F}.
\end{eqnarray*}
This algorithms requires $O(kmn)$ complexity (since the construction
of sampling probabilities depends on right singular vectors of $A$).
In \cite{mahoneydrineas}, the columns and rows were chosen randomly
according to a probability distribution formed by normalized statistical
leverage scores (based on right singular values of $A$). This algorithm
takes $A,$ $k$ and $\epsilon$ as input and uses column subset selection
procedure with $c=O(k\,logk\,\epsilon^{-2})$ columns of $A$ to construct
$C$ and with $r=O(k\,logk\,\epsilon^{-2})$ rows of $A$ to construct
$R.$ The matrix $U$ is given by $U=C^{\dagger}AR^{\dagger}.$ This
algorithm requires $O(kmn)$ complexity. In \cite{wangzhangli,wangzhang},
an improved algorithm has been proposed to obtain $CUR$ decomposition
with in shorter time compared to the existing relative error CUR algorithms
\cite{drineasmahoneymuthu,mahoneydrineas}.

The applicability of CUR decomposition in various fields can be found
in \cite{arisim,mittay,thukerst}. The generalization of CUR decomposition
to Tensors has been described in \cite{caifa}.\medskip{}

\subsubsection*{(II). Random projection based methods\medskip{}
}

The random projection method for low rank approximation of a matrix
is based on the idea of random projection. In random projection, the
original $d$ dimensional data is projected to a $k$ dimensional
($k<<d)$ subspace by post-multiplying a $k\times d$ random matrix
$\Omega$ (a matrix whose entries are independent random variables
of some specified distribution). The idea of random mapping is based
on the Johnson-Lindenstrauss lemma which says any set of $n$ points
in the $d$ dimensional Euclidean space can be embedded into $k$
dimensional Euclidean space such that the distance between the points
is approximately preserved \cite{dasgupata,johlinde}. 

The choice of the random matrix $\Omega$ plays an important role
in the random projection. There are several possible choices for $\Omega.$
The Bernoulli random matrix (with matrix entries 1 or -1 with an equal
probability of each), Gaussian random matrix (with matrix entries
have zero mean and unit variance normal distribution) are among the
choices for $\Omega.$ The details of several other choices for random
matrices were discussed in \cite{achlioptos,mahoney}.

The idea of the random projection based algorithms for low rank approximation
$\widetilde{A}$ of a matrix $A_{m\times n}$ is given below \cite{Halkotropp,mahoney,sapp}.
Let $k$ be the target rank, $s$ be the number of samples. 

Step 1. Consider a random matrix $\Omega_{n\times s}$.

Step 2. Obtain the product $Y_{m\times s}=A\Omega.$

Step 3. Compute an approximate orthonormal basis $Q_{m\times k}$
for the range of $Y$ via SVD.

Step 4. Finally obtain $\widetilde{A}=QQ^{T}A.$ 

In \cite{sarlos}, a structured random matrix has been considered
with $s=O(k/\epsilon)$ columns and low rank approximation has been
obtained such that
\begin{eqnarray*}
\left\Vert A-\widetilde{A}\right\Vert _{F}\leq(1+\epsilon)\,\left\Vert A-A_{k}\right\Vert _{F}
\end{eqnarray*}
 holds with high probability. The complexity of this algorithm is
$O(Mk/\epsilon+(m+n)k^{2}/\epsilon^{2}),$ where $M$ is the number
of non zero elements in $A$ and it requires 2 passes over the data.
In \cite{Halkotropp}, a standard Gaussian matrix has been considered
as $\Omega$ with $s=k+p$ columns, where $p\geq2$ an oversampling
parameter. The algorithm gives low rank approximation such that 
\begin{eqnarray*}
\left\Vert A-\widetilde{A}\right\Vert _{F}^{2}\leq(1+\frac{k}{p-1})\,\left\Vert A-A_{k}\right\Vert _{F}^{2}
\end{eqnarray*}
 holds with high probability. The complexity of this algorithm is
$O(mns+ms^{2}).$ 

This algorithm was further improved by coupling a form of the power
iteration method with random projection method \cite{Halkotropp,roklintygert}.
$Y$ in this modified algorithm is $Y=(AA^{T})^{q}A\Omega,$ where
$q$ is an iteration parameter. This provides the improved error estimates
of the form
\begin{eqnarray*}
\left\Vert A-\widetilde{A}\right\Vert _{F}^{2}\leq(1+\frac{k}{p-1})^{1/2q+1}\,\left\Vert A-A_{k}\right\Vert _{F}^{2}
\end{eqnarray*}
 with an extra computational effort. 

One can look at \cite{Halkotropp} for the error estimates in spectral
norm. The matrix multiplication $A\Omega$ requires $O(mns)$ operations
in the above algorithms. Some special structured matrices like $\Omega=DHS$
(details can be found in \cite{mahoney}) and subsampled random Fourier
transform (SRFT) matrices requires $O(mn\,logs)$ complexity \cite{Halkotropp}.
Complete analysis of the random projection methods can be found in
\cite{Halkotropp}. Random projection method also has been studied
in \cite{gu,ngundo,papavempala,roklintygert,vempala}.

\subsubsection*{(III). Randomized SVD}

The classical algorithms to compute SVD become very expensive as the
data size increases and also they require $O(k)$ passes over the
data. As explained in the introduction one can approximate SVD using
the randomized algorithms with less computational cost and fewer passes
over the data. It is a two stage procedure. In the first stage random
sampling is used to obtain a reduced matrix whose range approximates
the range of $A.$ The reduced matrix is factorized in the second
stage. This can be done in simple three steps \cite{Halkotropp}.

1. Form $B=Q^{T}A,$ which gives the low rank approximation $A_{k}=QB,$
where $Q=[q_{1},q_{2},...,q_{k}]$ is orthonormal basis obtained in
step 3 in the Random projection methods. 

2. Compute an SVD of the small matrix: $B=\tilde{U}\Sigma V^{T}.$

3. Set $U=Q\tilde{U}.$

This approximates the SVD with the same rank as the basis matrix $Q.$
The efficient implementation and approximation error of this procedure
can be found in section 5 of \cite{Halkotropp}. This scheme is well
suited for sparse and structured matrices. If the singular values
of $A$ decay slowly then power iterations (with $q=1$ or $2$) were
used (to form $Y$ in the random projection methods) to improve the
accuracy \cite{Halkotropp,roklintygert}. This gives the truncated
SVD (rank $k$) such that 
\begin{eqnarray*}
\left\Vert A-U\varSigma_{k}V^{T}\right\Vert _{2}\leq\sigma_{k+1}+\left[1+4\sqrt{\frac{2\,\textrm{min}\{m,n\}}{k-1}}\right]^{1/2q+1}\sigma_{k+1}
\end{eqnarray*}
holds with high probability. Here $k$ satisfies $2\leq k\leq0.5\,\textrm{min}\{m,n\}.$
The total cost of this algorithm to obtain rank $k$ SVD including
the operation count to obtain $Q$ is $O(mnlog(k)+k^{2}(m+n)).$ Randomized
SVD have also been studied and used in many applications \cite{dehdarideustch,drinearshuggins,jili,likwoklu2,martinsonszlamtyget,xiangzou,zhangerway}.$\blacksquare$ 

Some other randomized algorithms for low rank approximation of a matrix
have been proposed in \cite{Achliptossherry} (sparsification), \cite{bhojanapally,ubarumazum}.
Performance of different randomized algorithms have been compared
in \cite{bhojanapally,krishnamenon}.

\subsection{Some other techniques}

\subsubsection*{Non negative matrix factorization (NMF)}

Non negative matrix factorization of a given non negative matrix $A_{m\times n}$
(i.e all the matrix entries $a_{ij}\geq0)$ is finding two non negative
matrices $W_{m\times k}$ and $H_{k\times n}$ such that $WH$ approximates
$A.$ The chosen $k$ is much smaller than $m$ and $n.$ In general
it is not possible to obtain $W$ and $H$ such that $A=WH$. So NMF
is only an approximation. This problem can be stated formally as follows.\medskip{}
\\
\textbf{Definition (NMF problem): }Given a non negative matrix $A_{m\times n}$
and a positive integer $k<\textrm{min}\{m,n\},$ find non negative
matrices $W_{m\times k}$ and $H_{k\times n}$ to minimize the functional
\begin{eqnarray*}
f(W,H)=\frac{1}{2}\left\Vert A-WH\right\Vert _{F}^{2}.
\end{eqnarray*}

This is a nonlinear optimization problem. This factorization has several
applications in image processing, text mining, financial data, chemometric
and blind source separating etc. Generally, the factors $W$ and $H$
are naturally sparse, so they require very less storage. This factorization
has some disadvantages too. The optimization problem defined above
is convex in either $W$ or $H$, but not in both $W$ and $H$, which
means that the algorithms can only, if at all, guarantee the convergence
to a local minimum \cite{langvillemeyer}. The factorization is also
not unique (different algorithms gives different factorizations).

Such a factorization was first introduced in \cite{paaterotapper}
and the article \cite{leeseung1} about NMF became popular. There
are several algorithms available in the literature. Multiplicative
update algorithm \cite{leeseung1,leesung2}, projected gradient method
\cite{lin}, alternating least squares method \cite{cichradaf} and
several other algorithms described in \cite{berrybrowne},\cite{chikolda},\cite{gillis},\cite{kimparx},\cite{kimpark2}
and \cite{langvillemeyer} are among the algorithms for NMF. The non
negative tensor factorizations are described in \cite{flatz} and
several algorithms for both non negative matrix and tensor factorizations
with applications can be found in the book \cite{cichophan}.

\subsubsection*{Semidiscrete matrix decomposition (SDD)}

A semidiscrete decomposition (SDD) expresses a matrix as weighted
sum of outer products formed by vectors with entries constrained to
be in the set $\mathcal{S}=\{-1,0,1\}.$ The SDD approximation ($k$
term SDD approximation) of an $m\times n$ matrix $A$ is a decomposition
of the form \cite{koldaleary1} 
\begin{eqnarray*}
A_{k}=\underset{\underset{}{X_{k}}}{\underbrace{[x_{1}x_{2}...x_{k}]}}\underset{D_{k}}{\underbrace{\left[\begin{array}{cccc}
d_{1} & 0 & \cdots & 0\\
0 & d_{2} & \cdots & 0\\
\vdots &  & \ddots & \vdots\\
0 & \cdots & 0 & d_{k}
\end{array}\right]}}\underset{Y_{k}^{T}}{\underbrace{\left[\begin{array}{c}
y_{1}^{T}\\
y_{2}^{T}\\
\vdots\\
y_{k}^{T}
\end{array}\right]}}=\sum_{i=1}^{k}d_{i}x_{i}y_{i}^{T}.
\end{eqnarray*}
 Here each $x_{i}$ is an $m-$vector with entries from $\mathcal{S}=\{-1,0,1\},$
each $y_{i}$ is a $n-$vector with entries from the set $\mathcal{S}$
and each $d_{i}$ is a positive scalar. 

The columns of $X_{k},Y_{k}$ do not need to be linearly independent.
The columns can repeated multiple times. This $k$ term SDD approximation
requires very less storage compared to truncated SVD but it may require
large $k$ for accurate approximation. This approximation has applications
in image compression and data mining. This approximation was first
introduced in \cite{learypeleg} in the contest of image compression
and different algorithms have been proposed in \cite{koldaleary1,koldaleary2}.
A detailed description of SDD approximation with applications in data
mining can be found in the book \cite{skillicorn} and some other
applications can be found in \cite{luo,qiangxiao}.

\subsubsection*{Nystr$\textrm{\ensuremath{\ddot{o}}}$m Method}

The Nystr$\textrm{\ensuremath{\ddot{o}}}$m approximation is closely
related to $CUR$ approximation. Different from $CUR$, Nystr$\textrm{\ensuremath{\ddot{o}}}$m
methods are used for approximating the symmetric positive semidefinite
matrices (large kernel matrices arise in integral equations). The
Nystr$\textrm{\ensuremath{\ddot{o}}}$m method has been widely used
in machine learning community. The Nystr$\textrm{\ensuremath{\ddot{o}}}$m
method approximates the matrix only using a subset of its columns.
These columns are selected by different sampling techniques. The approximation
quality depends on the selection of the good columns. A brief description
of the Nystr$\textrm{\ensuremath{\ddot{o}}}$m approximation is given
below.

Let $A\in\mathbb{R}^{n\times n}$ be a symmetric positive semidefinite
matrix (SPSD). Let $C_{n\times m}$ be a matrix consists of $m$ ($<<n)$
randomly selected columns of $A$ as columns. Now the matrix $A$
can be rearranged such that $C$ and $A$ are written as 
\begin{eqnarray*}
C=\left[\begin{array}{c}
W\\
S
\end{array}\right]\,\,and\,\,\,A=\left[\begin{array}{cc}
W & S^{T}\\
S & B
\end{array}\right],
\end{eqnarray*}
 where $W\in\mathbb{R}^{m\times m},\,S\in\mathbb{R}^{(n-m)\times m}$
and $B\in\mathbb{R}^{(n-m)\times(n-m)}.$ Since $A$ is SPSD, $W$
is also a SPSD. 

For $k$ ($k\leq m$), the rank $k$ Nystr$\textrm{\ensuremath{\ddot{o}}}$m
approximation is defined by 
\begin{eqnarray*}
\widetilde{A}_{k}=CW_{k}^{\dagger}C^{T}
\end{eqnarray*}
 where $W_{k}$ is the best rank $k$ approximation of $W$ and $W_{k}^{\dagger}$
is the pseudoinverse of $W_{k}.$ $W_{k}^{\dagger}={\displaystyle \sum_{i=1}^{k}}\sigma_{i}^{-1}U^{(i)}(U^{(i)})^{T},$
where $\sigma_{i}$ is the $i^{th}$ singular value of $W$ and $U^{(i)}$
is the $i^{th}$ column of the matrix $U$ in the SVD of $W.$ The
computational complexity is $O(nmk+m^{3})$ which is much smaller
than the complexity $O(n^{3})$ of direct SVD. Using $W^{\dagger}$
instead of $W_{k}^{\dagger}$ gives more accurate approximation with
higher ranks than $k.$ 

The Nystr$\textrm{\ensuremath{\ddot{o}}}$m method was first introduced
in \cite{williseeger}. They have selected the random columns using
uniform sampling without replacement. A new algorithm has been proposed
and theoretically analyzed in \cite{drinieasmahoney_nystrom}. The
columns have been selected randomly using non-uniform probability
distribution and the error estimates of the Nystr$\textrm{\ensuremath{\ddot{o}}}$m
approximation were presented. By choosing $O(k/\epsilon^{4})$ columns
of $A,$ the authors have shown that 
\begin{eqnarray*}
\left\Vert A-CW_{k}^{\dagger}C^{T}\right\Vert _{\xi}\leq\left\Vert A-A_{k}\right\Vert _{\xi}+\epsilon\sum_{i=1}^{n}A_{ii}^{2}.
\end{eqnarray*}
 where $\xi=2,\,F$ and $A_{k}$ is the best rank $k$ approximation
of $A.$ 

These estimates have been further improved in \cite{gittensmahoney,talwalrostami,wangzhang,zhangtsang}.
A detailed comparison of the existing algorithms and error estimates
have been discussed in \cite{gittensmahoney}. The ensemble Nystr$\textrm{\ensuremath{\ddot{o}}}$m
method has been proposed in \cite{cortessanjivkumar,sanjeevmohri}.
Adaptive sampling techniques are used to select the random columns.
In \cite{likowklu}, a new algorithm which combines the randomized
low rank approximation techniques \cite{Halkotropp} and the Nystr$\textrm{\ensuremath{\ddot{o}}}$m
method was proposed. In \cite{Nemtsovaverbuch}, the details how the
Nystr$\textrm{\ensuremath{\ddot{o}}}$m method can be applied to find
the SVD of general matrices were shown.

\section{Cross/Skeleton approximation techniques}

In this section we discuss in detail the cross algorithms which gives
the low rank approximation of a matrix $A_{m\times n}.$ In these
algorithms the approximation of a matrix is obtained using the crosses
formed by the selected columns and rows of the given matrix. The computational
complexity of these algorithms is linear in $m,n$ and they use a
small portion of the original matrix. 

As described in the last section cross/skeleton approximation of a
matrix $A$ is given by $A\simeq CGR,$ where $C_{m\times k},R_{k\times n}$
consists of selected $k$ columns and $k$ rows of $A$ and $G=M^{-1},$
where $M_{k\times k}$ is the submatrix on the intersection of the
crosses formed by selected rows and columns from $A.$ In \cite{gorotyrzama1},
it has been shown that one can obtain a rank $k$ approximation within
an accuracy $\epsilon$ such that
\begin{eqnarray*}
\left\Vert A-CM^{-1}R\right\Vert _{2}=O\left(\epsilon\left\Vert A\right\Vert _{2}^{2}\left\Vert M^{-1}\right\Vert _{2}^{2}\right)
\end{eqnarray*}
 provided that $M$ is nonsingular. 

If $M$ is ill conditioned or if $M$ is singular then $CM^{-1}R$
will not approximate $A.$ So, the accuracy of the approximation depends
on the choice of $M.$ A good choice for $M$ is the maximum volume
submatrix i.e, the submatrix $M$ has determinant with maximum modulus
among all $k\times k$ submatrices of $A$ \cite{GoreTryty}. Since
the search for this submatrix is NP-complex problem \cite{civrilismail2},
it is not feasible even for moderate values of $m,n$ and $k.$ In
practice, such a submatrix $M$ can be replaced by matrices that can
be computed by the techniques, like adaptive cross approximation \cite{Bebendorf,BebendorfApp,bebendorf1},
skeleton decomposition with suboptimal maximum volume submatrix \cite{GoreOsele,GoreTryty}
and pseudoskeleton decomposition \cite{gorotyrzama1,Goritryzamar}.
The adaptive cross approximation (ACA) technique constructs the approximation
adaptively and the columns and rows are iteratively added until an
error criterion is reached. In pseudoskeleton decomposition the matrix
$G$ is not necessarily equal to $M^{-1}$ and even not necessarily
nonsingular. The applications of these techniques can be found in
\cite{BebendorfApp,bebendorf1,oseledets,Trytyshivkov}. Here we describe
all these techniques in detail.

\subsection{Skeleton decomposition}

\noindent Consider a matrix $A$ of order $m\times n.$ As described
above the matrix $A$ is approximated by $A\approx CGR,$ where $C$
and $R$ contains $k$ selected columns and rows respectively and
$G=M^{-1},$ where $M=A(I,J)$ is of order $k\times k$, a submatrix
on the intersection of the selected columns and rows ($I,\,J$ are
indices of rows and columns respectively). As explained above, obtaining
the maximum volume submatrix $M$ is very difficult, so we replace
it by a quasioptimal maximal volume submatrix. It has been discussed
in \cite{GoreOsele,GoreTryty}, how to find such a matrix. An algorithm
has been developed (named as \textquotedbl{}maxvol\textquotedbl{}
in \cite{GoreOsele}) and complexity of this algorithm has been shown
to be $O(mk).$ The algorithm takes a $m\times k$ matrix as input
and gives $k$ row indices such that the intersection matrix $M$
has almost maximal volume as output. To construct a rank $k$ skeleton
approximation of a given matrix $A$ using maxvol algorithm we follow
the steps given below. \medskip{}

\noindent \textit{Step 1}: Compute $k$ columns of $A$ given by the
indices $J=(j^{(1)},j^{(2)},\ldots,j^{(k)})$ and store them in a
matrix. Let it be $C$. $i.e.,C=A(:,J)$ is of order $m\times k.$

\noindent \textit{Step 2}: Now we find a good matrix by using the
maxvol procedure on $C$ \cite{GoreOsele} that gives $k$ row indices
$I=(i^{(1)},i^{(2)},\ldots,i^{(k)})$ such that the corresponding
intersection matrix, say $M$ has almost maximal volume.

\noindent \textit{Step 3}: Store the $k$ rows in a matrix $R=A(I,:)$
($R$ is of order $k\times n$).

\noindent \textit{Step 4}: Therefore the skeleton decomposition is
$A\approx CGR$ where $G=M^{-1},$ $M=A(I,J)$ is of order $k\times k.$\medskip{}

\noindent \textbf{\textit{$Remark$}}\textbf{: }\textit{ If the column
indices $J=(j^{(1)},j^{(2)},\ldots,j^{(k)})$ at the beginning were
badly chosen (may be because of some random strategy), this approximation
might not be very good and the inverse of $M$ might be unstable (nearly
singular matrix). Even if the ranks are over estimated the inverse
of $M$ will also be unstable \cite{oseledets}.}\medskip{}

\noindent To overcome this, after getting good row indices $I,$ one
can use the maxvol procedure for the row matrix $R$ to optimize the
choice of the columns and even alternate further until the determinant
of $M$ stays almost constant and the approximation is fine \cite{oseledets,Trytyshivkov}.
In the case of over estimation of ranks, some remedy is still possible
to get good accuracy \cite{oseledets}.

\medskip{}
\textbf{Computational complexity: }The complexity of the algorithm
is $O\left((m+n)k^{2}\right),$ and only $k(m+n)$ of the original
entries of $A$ have to be computed. The storage required for the
approximation is $k(m+n).$ We can see from the algorithm that only
few original matrix entries have been used in the final approximation.
\medskip{}

A quasioptimal error estimate for skeleton approximation of matrix
has been derived in \cite{gorityr3}. It has been shown that if the
matrix $M$ has maximal in modulus determinant among all $k$ by $k$
submatrices of $A$ then 
\begin{eqnarray*}
\left\Vert A-CM^{-1}R\right\Vert _{\infty}\leq(k+1)^{2}\begin{array}{c}
\textrm{min}\\
\textrm{rank}\,B\leq k
\end{array}\left\Vert A-B\right\Vert _{\infty}.
\end{eqnarray*}
 Where $\left\Vert A\right\Vert _{\infty}$ defined as the largest
entry in the absolute value of the matrix $A$ (sup-norm).

\subsection{Pseudoskeleton approximation}

\noindent As explained earlier, the intersection matrix $M$ need
not be invertible in the pseudoskeleton approximation. The theory
of pseudoskeleton approximation has been studied in \cite{gorotyrzama1,gorotyrzamar,Goritryzamar}.
For example, $G$ can be chosen as the pseudoinverse of $M.$ Here
we describe a similar algorithm presented in \cite{espigkumarjan},
which can be used for the treatment of higher dimensional arrays.
To construct pseudoskeleton approximation of $A_{m\times n}$ with
rank $r\leq\textrm{min}(m,n)$, we proceed as follows.

\noindent Instead of choosing $r$ columns from the matrix $A,$ we
choose $k$ ($k>r)$ random columns given by the indices $J=(j^{(1)},j^{(2)},\ldots,j^{(k)})$
from the given matrix $A$ and store them in a matrix, let it be $C,$
therefore $C=A(:,J)\in\mathbb{R}^{m\times k}.$ Now we use maximum
volume submatrix procedure on $C$ \cite{GoreOsele} to find $k$
row indices $I=(i^{(1)},i^{(2)},\ldots,i^{(k)})$ corresponding to
$C$. Store the $k$ row indices in a matrix $R.$ Therefore $R=A(I,:)\in\mathbb{R}^{k\times n}$
and the intersection matrix $M=A(I,J)$ of $C$ and $R$ has almost
maximal volume. Now the approximation looks like 
\[
A\sim C_{m\times k}M_{k\times k}^{-1}R_{k\times n}.
\]
Due to the random choice, inverse of $M$ may not be stable. To overcome
this, now we decompose the matrix $M$ using the Singular Value Decomposition.
i.e., 
\[
M=U_{M}S_{M}V_{M}^{T}.
\]

\noindent Where $U_{M},$ $V_{M}$ are orthogonal matrices and $S_{M}$
is a diagonal matrix contains $k$ singular values $\sigma_{1}\geq\sigma_{2}\geq\cdots\geq\sigma_{k}.$
Now we find singular values which are greater than $\epsilon$ (this
value is predefined and similar to $\tau$ used in \cite{gorotyrzamar})\footnote{Different conditions can be considered.}
from diagonal matrix $S_{M}.$ Let $r$ be the number of singular
values larger than $\epsilon$ and replace the other singular values
by zero. Truncate the matrices $U_{M},$ $S_{M}$ and $V_{M}$ according
to the $r$ singular values and store them as $U_{r}\in\mathbb{R}^{k\times r},$
$V_{r}\in\mathbb{R}^{k\times r}$ and $S_{r}\in\mathbb{R}^{r\times r}.$
Therefore the pseudoinverse of $M$ is 
\[
M^{\dagger}=(U_{r}S_{r}V_{r}^{T})^{-1}=V_{r}S_{r}^{-1}U_{r}^{T}.
\]

\noindent Therefore the matrix $A$ is approximated as 
\[
A\sim CV_{r}S_{r}^{-1}U_{r}^{T}R.
\]

\noindent So the approximation of $A$ looks as 
\[
A\sim CR
\]

\noindent where $C=CV_{r}S_{r}^{-1},$ $R=U_{r}^{T}R.$

One can further compress the rank $r$ using the technique described
in \cite{BebendorfApp}. The above described procedure is similar
to the method proposed in \cite{ZhuXlinW} where both the columns
and rows are chosen randomly.

\subsection*{Algorithm\medskip{}
}

\noindent \textit{Step 1}: Choose $k$ columns randomly from the matrix
$A$ given by the indices

\noindent $J=(j^{(1)},j^{(2)},\ldots,j^{(k)})$ and store them in
a matrix. Let it be $C$. $i.e.,C=A(:,J)$ is of order $m\times k.$

\noindent \textit{Step 2}: Now use maxvol procedure on $C$ \cite{GoreOsele}
to get $k$ row indices $I=(i^{(1)},i^{(2)},\ldots,i^{(k)})$ from
$A$ corresponding to the columns in $C$ and store them in a matrix,
say $R=A(I,:)$ such that the corresponding intersection matrix, say
$M=A(I,J)$ has almost maximal volume.

\noindent \textit{Step 3}: Apply SVD to the matrix $M$ and let $U_{M},$
$V_{M}$ and $S_{M}$ be the decomposition matrices by SVD of $M.$

\noindent \textit{Step 4}: Fix $\epsilon$ and find $r$ singular
values satisfying the condition $\sigma_{i}>\epsilon$ for $i=1,2,\ldots,k.$

\noindent \textit{Step 5}: Now truncate the matrices $U_{M},$ $V_{M}$
and $S_{M}$ according to $r$ singular values and store them in the
matrices $U_{r},$ $V_{r}$ and $S_{r}$ respectively.

\noindent \textit{Step 6}: Find the pseudoinverse of $M.$ i.e., $M^{\dagger}=(U_{r}S_{r}V_{r}^{T})^{-1}=V_{r}S_{r}^{-1}U_{r}^{T}.$

\noindent \textit{Step 7}: So, finally $A$ will be decomposed as
$A\sim CR,$ where $C=CV_{r}S_{r}^{-1},$ $R=U_{r}^{T}R.$\bigskip{}

\noindent \textbf{Computational Complexity:} The overall computational
cost of the above algorithm is

\noindent 
\[
O(k^{3}+rk(m+n)+mr^{2}).
\]

Here one can see that the most dominating one in operation count is
$O(k^{3})$ (due to the SVD on $M$ in step 3). Some terms in the
complexity involves $mrk,\,nrk$ (due to the computations in step
7). Since $k\ll m,n$ and $r\ll k$ they do not dominate. So the overall
complexity of the algorithm is linear in $m$ and $n.$

Pseudoskeleton approximation is used to construct different tensor
decomposition formats \cite{espigkumarjan,oseledets}. Different applications
of pseudoskeleton approximation can be found in \cite{carin,ZhuXlinW}.\medskip{}
\\
\textbf{Error Analysis: }In \cite{gorotyrzama1,Goritryzamar}, the
error in the pseudoskeleton approximation has been studied. It has
been shown that if the matrix $A$ is approximated by rank $r$ within
an accuracy $\epsilon$ then there exists a choice of $r$ columns
and $r$ rows i.e $C$ and $R$ and the intersection matrix $M$ such
that $A\simeq CGR$ satisfying 
\begin{eqnarray*}
\left\Vert A-CGR\right\Vert _{2}\leq\epsilon\left(1+2\sqrt{mr}+2\sqrt{nr}\right).
\end{eqnarray*}

Here the columns of $C$ and rows of $R$ are chosen such that their
intersection $M$ has maximal volume \cite{GoreTryty}. 

In \cite{chiulaurent} sublinear randomized algorithm for skeleton
decomposition has been proposed. Uniformly sampling $l\simeq r\,log\,(\textrm{max}(m,n))$
rows and columns are considered to construct a rank $r$ skeleton
approximation. The computational complexity of the algorithm is shown
to be $O(l^{3})$ and the following error estimate has been proved. 

Suppose $A\simeq X_{1}A_{11}Y_{1}^{T}$ where $X_{1},Y_{1}$ have
$r$ orthonormal columns (not necessarily singular vectors of $A$)
and $A_{11}$ is not necessarily diagonal. Assume $X_{1}$ and $Y_{1}$
are incoherent, then 
\begin{eqnarray*}
\left\Vert A-CGR\right\Vert _{2}=O\left(\frac{\sqrt{nm}}{l}\left\Vert A-X_{1}A_{11}Y_{1}^{T}\right\Vert _{2}\right)
\end{eqnarray*}
holds with high probability.

\subsection{Adaptive Cross Approximation (ACA)}

Adaptive cross approximation has been introduced in \cite{Bebendorf,bebendorf1}.
In contrast to the pseudoskeleton method here the rows and columns
are chosen adaptively such that in each step a rank one approximation
is added to the approximant. We try to keep the notation used when
it was first investigated thoroughly, so we prefer the language of
functions here instead of matrices identifying both by a regular,
sufficiently dense grid in $[0,1]^{2}$, i.e. 
\[
A=(a_{ij}),\quad\text{ where }\quad a_{ij}=f\left(\frac{i-1}{n-1},\frac{j-1}{n-1}\right)
\]
for indices $1\leq i,j\leq n$ and a given function $f:[0,1]^{2}\rightarrow\mathbb{R}$.\\
 We are basically concerned with the questions: How to approximate
$f$ by something like 
\begin{equation}
f\sim\sum_{i=1}^{k}g_{i}\otimes h_{i},\label{basis}
\end{equation}
i.e., by a finite sum of tensor products of one-dimensional functions
(here we write ($g\otimes h)(x,y)=g(x)h(y)$)? And how good is this
approximation?\\
 The first famous result in this direction is due to Schmidt \cite{schmidt},
who gave a complete answer in the case $f\in L_{2}$. A standard reference
for questions in this area is \cite{cheney2}, a nice survey can be
found in \cite{cheney1}.\\
 Now a very special choice of functions $g,h$ in \eqref{basis} is
considered, namely the restriction of $f$ itself to certain lines.
In the discrete setting, that means we only allow columns and rows
of the matrix $A$ to be building blocks.

Let $f:[0,1]^{2}\longrightarrow\mathbb{R}$, then the recursion $R_{0}(x,y)=f(x,y)$
and 
\begin{equation}
R_{k}(x,y)=R_{k-1}(x,y)-\frac{R_{k-1}(x,y_{k})R_{k-1}(x_{k},y)}{R_{k-1}(x_{k},y_{k})}\quad\text{for }k\in\mathbb{N},\label{heart}
\end{equation}
with points $1\leq x_{k},y_{k}\leq n$ chosen such that $R_{k-1}(x_{k},y_{k})\neq0$,
is the heart of the two-dimensional cross approximation, compare \cite{bebendorf1,Schneider}.
So in each step a pivot $(x_{k},y_{k})$ with $f(x_{k},y_{k})\neq0$
is chosen (the process of choosing these pivots is called pivoting)
and the corresponding row and column is used to add another rank $1$
approximation for the remainder $R_{k-1}$. After $k$ iteration steps
the approximant $S_{k}(x,y)=f(x,y)-R_{k}(x,y)$ is calculated. $S_{k}(x,y)$
takes the form \cite{Bebendorf,bebendorf1}
\begin{eqnarray*}
S_{k}(x,y)=\sum_{i,j=1}^{k}(M_{k})_{ij}^{-1}f(x,y_{i})f(x_{j},y)
\end{eqnarray*}
 where $(M_{k})_{ij}=f(x_{i},y_{j}),\,\,\,i,j=1,2,..,k.$ 

In the matrix form, the columns of $C$ and rows of $R$ are iteratively
added and the approximation of $A_{n\times n}$ takes the form (compare
with pseudoskeleton approximation) 
\begin{eqnarray*}
A\simeq CGR,\,\,\textrm{where} &  & G=M_{k}^{-1}.
\end{eqnarray*}

The cross approximation has nice properties like interpolation property
and rank property. For any function $f:[0,1]^{2}\rightarrow\mathbb{R}$,
we have the interpolation property 
\[
R_{k}(x,y)=0,\quad\text{ as long as }\quad x=x_{i}\vee y=y_{i}\quad\text{ for at least one }\quad i\in\{1,\ldots,n\}.
\]
That means, on the chosen crosses this procedure is exact\textbf{.}
The next result takes an a priori knowledge about structural properties
of the underlying function into account. We say that a function $f$
has separation rank $k$, if one can represent it as 
\[
f(x,y)=\sum_{i=1}^{k}g_{i}(x)h_{i}(y)
\]
and there is no such representation with reduced summing order. This
is just the continuous analog to the rank of a matrix. We call the
following the rank property: If $f$ has separation rank $k$ cross
approximation reproduces $f$ after $k$ steps exactly, that means
\[
R_{k}=f-S_{k}=0\quad\text{on }[0,1]^{2}.
\]
A matrix version of this result was first proved in \cite{bebendorf1}(Lemma
7).\\
 There also has been some effort to the error analysis: In \cite{Schneider}
was stated 
\begin{equation}
|R_{k}(x,y)|\leq(k+1)^{2}E(f,G)_{C([0,1]^{2})},\label{error}
\end{equation}
where $E(f,G)_{C([0,1]^{2})}=\inf_{g\in G}\|f-g\|_{\infty}$ is the
error of best approximation of $f$ in 
\[
G=\left\{ g=\sum_{i=1}^{k}\varphi_{i}(x)\psi_{i}(y),\quad\varphi_{i},\psi_{i}\in C([0,1])\right\} 
\]
measured in the $\sup$-norm (compare with the matrix version of this
result stated at the end of section 3.1). Similar results can be found
in \cite{Bebendorf,BebendorfApp}. In \eqref{error} a very special
choice of pivots is crucial, namely the maximal volume concept, i.e.
$(x_{1},y_{1}),\ldots,(x_{k},y_{k})$ are chosen, such that 
\[
\left|\det\left(f(x_{i},y_{j})\right)_{i,j=1}^{k}\right|
\]
is maximal under all possible choices of the points. This is of course
not practical and since one wants to keep the pivots of the previous
steps untouched, a good alternative is the partial pivoting. Here
the positions in one direction are chosen by some strategy (for example
completely random) and in the second direction the maximum in modulus
of the remainder is taken on a line. This is still inexpensive and
leads to good numerical results, see \cite{Bebendorf,BebendorfApp,bebendorf1,naraparajujan}.

Practically it is expensive to update the whole remainder at each
iteration step. As described in \cite{Bebendorf,bebendorf1} the approximation
of the form $S_{k}={\displaystyle \sum_{i=1}^{k}\alpha_{i}u_{i}(x)v_{i}(y)}$
can be obtained without updating the whole remainder. 

If we set $u_{k}(x)=R_{k-1}(x,y_{k})$ and $v_{k}(y)=R_{k-1}(x_{k},y)$,
then using \eqref{heart} we get 
\begin{eqnarray}
S_{k}(x,y) & = & f(x,y)-R_{k}(x,y)\nonumber \\
 & = & f(x,y)-\left(R_{k-1}(x,y)-\frac{u_{k}(x)}{u_{k}(x_{k})}v_{k}(y)\right)\nonumber \\
 & \vdots\nonumber \\
 & = & f(x,y)-\left(R_{0}(x,y)-\frac{u_{1}(x)}{u_{1}(x_{1})}v_{1}(y)-\cdots-\frac{u_{k}(x)}{u_{k}(x_{k})}v_{k}(y)\right)\nonumber \\
 & = & \sum_{i=1}^{k}\frac{u_{i}(x)}{u_{i}(x_{i})}v_{i}(y),\label{2approx}
\end{eqnarray}
where we realize that it has the desired tensor product structure
of separated variables (compare with (16)). By a similar calculation
one can even derive the explicit formulas 
\begin{equation}
u_{k}(x)=f(x,y_{k})-\sum_{l=1}^{k-1}\frac{u_{l}(x)}{u_{l}(x_{l})}v_{l}(y_{k})\label{explu}
\end{equation}
and 
\begin{equation}
v_{k}(y)=f(x_{k},y)-\sum_{l=1}^{k-1}\frac{u_{l}(x_{k})}{u_{l}(x_{l})}v_{l}(y),\label{explv}
\end{equation}
see also \cite{Bebendorf}.\\
 For matrices that means instead of storing $n^{2}$ values one needs
to call the function less than $2kn$ times and store this reduced
amount of data. \\

\smallskip{}

Here we describe the construction of rank-$k$ approximation of the
matrix $A=(a_{ij})$ using ACA (described in \cite{Bebendorf,bebendorf1}).
\smallskip{}
\\
\textbf{Notation}: In the equations (20) and (21), the $u_{k}(x)$
and $v_{k}(y)$ are functions of $x$ and $y$ respectively. In the
matrix case, $\mathbf{u}_{k}=(u_{k}(1),....u_{k}(n))^{T}$ represents
a column vector and $\mathbf{v}_{k}=(v_{k}(1),....,v_{k}(n))$ denotes
a row vector. \medskip{}

\noindent \textbf{Details and operation count of the algorithm\medskip{}
}

\noindent We choose $y_{1}$ randomly at the beginning. In the first
step we choose the column vector $\mathbf{u}_{1}$ with entries $u_{1}(i)=A(i,y_{1}),i=1,2,..,n$
and find the maximum element index from $\left|\mathbf{u}_{1}\right|.$
\footnote{$\left|\mathbf{u}\right|=(\left|u(1)\right|,....,\left|u(n)\right|)$}
Let the maximum element index from $\left|\mathbf{u}_{1}\right|$
be $x_{1}$ and let $\delta_{1}=u_{1}(x_{1})$. Now find the corresponding
row vector $\mathbf{v}_{1}$ with entries $v_{1}(i)=A(x_{1},i),i=1,2,..,n.$ 

\noindent Now we find the maximum element index from $\left|\mathbf{v}_{1}\right|$
(the index should not be $y_{1}$) and let it be $y_{2}.$ In the
second step (here we know the pivot element $y_{2}$), we find the
vector $\mathbf{u}_{2}$ with entries (see equation 20)

\noindent 
\[
u_{2}(i)=A(i,y_{2})-\left((u_{1}(i)v_{1}(y_{2})\right)/\delta_{1}\,\,i=1,2,...,n.
\]
Let the maximum element index in $\left|\mathbf{u}_{2}\right|$ (the
index should not be $x_{1}$) be $x_{2}$ and let $\delta_{2}=u_{2}(x_{2}).$
Now it is easy to find the vector $\mathbf{v}_{2}$ with (see equation
21)

\noindent 
\[
v_{2}(i)=A(x_{2},i)-\left(u_{1}(x_{2})v_{1}(i)\right)/\delta_{1},\,\,i=1,2,...,n.
\]
So, here we are doing $n$ multiplications, $n$ subtractions and
$1$ division for $u_{2}$ and $v_{2}$ respectively. The total number
of operations at second step are $4n+2$. 

Let $y_{3}$ \footnote{In general the index at $k$th step should not be any index obtained
at previous steps } be the maximum element index in $\left|\mathbf{v}_{2}.\right|$ In
the third step, we find the vector $\mathbf{u}_{3}$ corresponding
to known $y_{3}.$ Therefore the vector 

\noindent 
\[
u_{3}(i)=A(i,y_{3})-\left((u_{1}(i)v_{1}(y_{3})\right)/\delta_{1}-\left(u_{2}(i)v_{2}(y_{3})\right)/\delta_{2},i=1,2,..,n.
\]
Let the maximum element index from $\left|\mathbf{u}_{3}\right|$
be $x_{3}$ and denote $\delta_{3}=u_{3}(x_{3}).$ Then the entries
of $\mathbf{v}_{3}$ are given by 

\noindent 
\[
v_{3}(i)=A(x_{3},i)-\left((u_{1}(x_{3})v_{1}(i)\right)/\delta_{1}-\left(u_{2}(x_{3})v_{2}(i)\right)/\delta_{2},\,i=1,2,..,n.
\]
So, here we are doing $2n$ multiplications, $n$ additions, $n$
subtractions and $2$ divisions for $\mathbf{u}_{3}$ and $\mathbf{v}_{3}$
respectively. The total number of operations at this step are $8n+4$.

\noindent Similarly we can find for other chosen rows/columns. Therefore
at the $k$'th step, we do $(k-1)(4n+2)$ operations.

\noindent The total number of operations in the algorithm is 

\noindent 
\[
0+4n+2+8n+4+\cdots+(k-1)(4n+2)=(4n+2)(1+2+\cdots+k-1)
\]

\noindent 
\[
\,\,\,\,\,\,\,\,\,\,\,\,\,\,\,\,\,\,\,\,\,\,\,\,\,\,\,\,\,\,\,\,\,\,\,\,\,\,\,\,\,\,\,\,\,\,\,\,\,\,\,\,\,\,\,\,\,\,\,\,\,\,\,\,\,\,\,\,\,\,\,\,\,\,\,\,\,\,\,\,\,\,\,=(4n+2)\frac{(k-1)k}{2}
\]

\noindent 
\[
\,\,\,\,\,\,\,\,\,\,\,\,\,\,\,\,\,\,\,\,\,\,\,\,\,\,\,\,\,\,\,\,\,\,\,\,\,\,\,\,\,\,\,\,\,\,\,\,\,\,\,\,\,\,\,\,\,\,\,\,\,\,\,\,\,\,\,\,\,\,\,\,\,\,\,\,\,\,\,\,\,=(2n+1)(k^{2}-k)
\]

\noindent 
\[
\,\,\,\,\,\,\,\,\,\,\,\,\,\,\,\,\,\,\,\,\,\,\,\,\,\,\,\,\,\,\,\,\,\,\,\,\,\,\,\,\,\,\,\,\,\,\,\,\,\,\,\,\,\,\,\,\,\,\,\,\,\,\,\,\,\,\,\,\,\,\,\,\,\,\,\,\,\,\,\,\,\,\,\,\,\,\,\,\,\,\,\,\,=2nk^{2}+k^{2}-2kn+k.
\]

\noindent Therefore, the complexity of this algorithm is $O(k^{2}n),$
which is linear in $n\,(k\ll n)$. In particular, if a matrix is of
dimension $m\times n,$ the number of operations required to construct
${\displaystyle S_{k}=\sum_{i=1}^{k}\frac{1}{\delta_{k}}\mathbf{u}_{k}\mathbf{v}_{k}}$
is $O(k^{2}(m+n)),$ while the storage required for the approximation
$S_{k}$ is of order $k(m+n)$.

An algorithm called Cross-2D has been described in \cite{oselsavostyr}.
The extension of adaptive cross approximation to higher order tensors
has been discussed in \cite{Bebendorf,naraparajujan} and different
applications of ACA can be found in \cite{aminfarambikadarve,BebendorfApp,bebendorf1,bebenrja,fladkhrom,hackbusch,kurzrain,rjasanowstien,rogus,zhaovouvlee}. 

We conclude this section with a remark on the performance of cross/skeleton
approximation techniques on the matrices with non smooth data. These
matrices generally arises from the discretization of the singular
functions in some applications \cite{Bebendorf,bebendorf1,oselsavostyr}.
The cross/skeleton approximation algorithms require larger ranks to
approximate the matrices with non smooth data \cite{Bebendorf,oselsavostyr}
compared to the ranks required to approximate the matrices with smooth
data.

\subsubsection*{Acknowledgements}

The authors wishes to thank Prof. Wolfgang Hackbusch for his constant
support and encouragement during their stay at Max-Planck institute
for Mathematics in the Sciences, Leipzig, Germany and National Board
of Higher Mathematics, India for financial support.

\end{document}